\documentclass{amsart}

\usepackage{amssymb,amsmath,color}
\usepackage{amsthm}
\usepackage{url}
\usepackage{tikz}

\definecolor{red}{rgb}{1,0,0}
\definecolor{blue}{rgb}{.2,.2,.8}

\def\ll{\lambda}
\def\m{\mu}
\def\T{\mathcal T_r}

\def\P{\mathcal P}
\def\C{\mathcal D_{1,r}}
\def\A{\mathcal O_{1,r}}
\def\D{\mathcal D_r}
\def\OO{\mathcal O_r} 
\def\MD{\mathcal{MD}}
\def\MDS{\mathcal{DD}_r}
\def\OD{\overline{\mathcal{D}}_r}
\def\ODT{\overline{\mathcal{D}}_3}
\def\Z{\mathbb Z}
\def\N{\mathbb N}

\newtheorem{theorem}{Theorem}[section]

\newtheorem{proposition}{Proposition}

\theoremstyle{definition}
\newtheorem{definition}{Definition}
\newtheorem{example}{Example}
\newtheorem{remark}{Remark}

\newcommand{\ds}{\displaystyle}

\begin{document}

\title{Beck-type identities for Euler pairs of order $r$}
\author{Cristina Ballantine}\address{Department of Mathematics and Computer Science\\ College of the Holy Cross \\ Worcester, MA 01610, USA \\} 
\email{cballant@holycross.edu} 
\author{Amanda Welch} \address{Department of Mathematics and Computer Science\\ College of the Holy Cross \\ Worcester, MA 01610, USA \\} \email{awelch@holycross.edu}

\maketitle


\begin{abstract} 
Partition identities are often statements asserting that the set $\P_X$ of partitions of $n$ subject to condition $X$ is equinumerous to the set $\P_Y$ of partitions of $n$ subject to condition $Y$. A Beck-type identity is a companion  identity to $|\P_X|=|\P_Y|$ asserting that the difference $b(n)$ between the number of parts in all partitions in $\P_X$ and the number of parts in all partitions in $\P_Y$ equals a $c|\P_{X'}|$ and also $c|\P_{Y'}|$, where $c$ is some constant related to the original identity, and $X'$, respectively $Y'$, is a condition on partitions that is a very slight relaxation of  condition $X$, respectively  $Y$. A second Beck-type identity involves the difference $b'(n)$ between the total number of different parts in all partitions in $\P_X$ and the total number of different parts in all partitions in $\P_Y$. We extend these results to Beck-type identities accompanying all identities given by Euler pairs of order $r$ (for any $r\geq 2$). As a consequence, we obtain many families of new Beck-type identities. We give analytic and bijective proofs of our results. 

\end{abstract}

{\bf Keywords:} partitions, Euler pairs, Beck-type identities, words

{\bf MSC 2010:}  05A17, 11P81, 11P83

\section{Introduction}

The origin of this article is rooted in  two conjectures by  Beck which appeared in The On-Line Encyclopedia of Integer Sequences \cite{B1} on the pages for sequences A090867 and A265251. The conjectures, as formulated by Beck, were proved by Andrews in \cite{A17} using generating functions. Certain generalizations and  combinatorial proofs appeared in \cite{FT17} and \cite{Y18}. Combinatorial proofs of the original conjectures were also given in \cite{BB19}. Several additional similar identities were proved in the last two years. See for example \cite{AB19, LW19, LW19b,LW19c}. In order to define Beck-type identities, we first introduce the necessary terminology and notation.

In this article $\N$ denotes the set of positive integers. Given a non-negative integer $n$, a \textit{partition} $\lambda$ of $n$ is a non-increasing sequence of positive integers $\lambda=(\lambda_1, \lambda_2, \ldots, \lambda_k)$ that add up to $n$, i.e., $\ds\sum_{i=1}^k\lambda_i=n$. Thus, if $\ll=(\lambda_1, \lambda_2, \ldots, \lambda_k)$ is a partition, we have $
\ll_1\geq \ll_2\geq \ldots \geq \ll_k$. The numbers $\lambda_i$ are called the \textit{parts} of $\lambda$ and $n$ is called the \textit{size} of $\lambda$. The number of parts of the partition is called the \textit{length} of $\lambda$ and is denoted by $\ell(\lambda)$. 

If $\ll, \m$ are two arbitrary partitions, we denote by $\ll\cup \m$ the partition obtained by taking all  parts  of $\ll$ and all parts of $\m$ and rearranging them to form a partition. For example, if $\ll=(5, 5, 3, 2, 2, 1)$ and $\m=(7, 5, 3, 3)$, then $\ll\cup \m=(7,5,5,5,3,3,3,2,2,1)$. 

When convenient,  we use the exponential notation for parts in a partition. The exponent of a part is the multiplicity of the part in the partition. For example, $(7, 5^2, 4, 3^3, 1^2)$ denotes the partition $(7, 5,5, 4, 3, 3, 3, 1, 1)$. It will be clear from the context when exponents refer to multiplicities and when they are exponents in the usual sense.

\medskip

Let $S_1$ and $S_2$ be subsets of the positive integers. We define 
$\OO(n)$ to be  the set of partitions of $n$ with all parts from the set $S_2$ and  $\D(n)$ to be the set of partitions of $n$ with  parts in $S_1$ repeated at most $r-1$ times.  Subbarao \cite{S71} proved the following theorem. 

\begin{theorem} \label{sub-thm} $|\OO(n)|=|\D(n)|$ for all non-negative integers $n$ if and only if $rS_1\subseteq S_1$ and $S_2=S_1\setminus rS_1$. \end{theorem} Andrews \cite{A69}  first discovered this result for $r=2$ and called a pair $(S_1,S_2)$ such that $|\mathcal{O}_2(n)|=|\mathcal{D}_2(n)|$  an \textit{Euler pair} since the pair $S_1=\N$ and $S_2=2\N-1$ gives  Euler's identity. By analogy, Subbarao called a pair $(S_1,S_2)$ such that $|\OO(n)|=|\D(n)|$  an \textit{Euler pair of order $r$}.

\begin{example}[Subbarao \cite{S71}]\label{eg:epair}
Let 
\begin{align*}
S_1 &= \{ m \in \N : \ m \equiv 1 \ (mod \ 2) \};\\
S_2 &= \{ m \in \N : \ m \equiv \pm 1 \ (mod \ 6) \}.
\end{align*}

\noindent Then $(S_1, S_2)$ is an Euler pair of order 3.
\end{example}

Note that Glaisher's bijection  used to prove $|\OO(n)|=|\D(n)|$ when $S_1=\N$ and $S_2=2\N-1$ can be generalized to any Euler pair of order $r$. If $(S_1,S_2)$ is an Euler pair of order $r$,  let $\varphi_r$ be the map from $\OO(n)$ to $\D(n)$ which repeatedly merges $r$ equal parts into a single part until there are no parts repeated more than $r-1$ times. The map $\varphi_r$ is a bijection and we refer to it as Glaisher's bijection. 

Given $(S_1,S_2)$, an Euler pair of order $r$, we refer to the elements in $S_2=S_1\setminus rS_1$ as \textit{primitive} elements and to the elements of $rS_1=S_1\setminus S_2$ as \textit{non-primitive} elements. We usually denote primitive parts by bold lower case letters, for example $\bf a$. Non-primitive parts are denoted by (non-bold) lower case letters. If $a$ is a non-primitive part of a partition and we want to emphasize the largest power $k$ of $r$ such that $a/r^k\in S_1$, we write $a=r^k\bf a$ with ${\bf a}\in S_2$ and $k \geq 1$. 

 Let $\A(n)$ be the set of partitions of $n$ with parts in $S_1$ such that the set of parts in $rS_1$ has exactly one element. Thus, a partition in $\A(n)$ has exactly one  non-primitive part (possibly repeated).   Let $\C(n)$ be the set of partitions of $n$ with  parts in $S_1$ in which exactly one part is repeated at least $r$ times. 

Let $a_r(n)=|\A(n)|$ and $c_r(n)=|\C(n)|$. Let $b_r(n)$ be the difference between the number of parts in all partitions in $\OO(n)$ and the number of parts in all partitions in $\D(n)$.  Thus, $$b_r(n)=\sum_{\lambda\in \OO(n)} \ell(\lambda)-\sum_{\lambda\in \D(n)} \ell(\lambda).$$

Let $\T(n)$ be the subset of $\C(n)$ consisting of partitions of $n$ in which  one part is repeated more than $r$ times but less than $2r$ times. Let $c'_r(n)=|\T(n)|$. Let $b'_r(n)$ be the difference between the total number of different parts  in all partitions in $\D(n)$ and the total number of different parts in all partitions in $\OO(n)$ (i.e., in each partition, parts are counted without multiplicity). If we denote by $\bar\ell(\ll)$ the number of different parts in $\ll$, then $$b'_r(n)=\sum_{\lambda\in \D(n)} \bar\ell(\lambda)-\sum_{\lambda\in \OO(n)} \bar\ell(\lambda).$$

In \cite{B1}, Beck conjectured that, if $S_1=\N$ and $S_2=2\N-1$, then $$a_2(n)=b_2(n)=c_2(n)$$ and $$c'_2(n)=b'_2(n).$$ Andrews proved these identities in \cite{A17}  using generating functions.  Combinatorial proofs were given in \cite{BB19}.  For the case $r\geq 2$,  $S_1=\N$, and  $S_2=\{k\in \N: k \not \equiv 0 \pmod r\}$, Fu and Tang \cite{FT17} gave generating function proofs for \begin{equation}\label{beck1} a_r(n)=\frac{1}{r-1}b_r(n)=c_r(n)\end{equation} and \begin{equation} \label{beck2} c'_r(n)=b'_r(n).\end{equation} They also proved combinatorially that  $a_r(n)=c_r(n)$.  In \cite{Y18},  Yang gave combinatorial proofs of \eqref{beck1} and \eqref{beck2} in the case $r\geq 2$,  $S_1=\N$, and  $S_2=\{k\in \N:  k \not \equiv 0 \pmod r\}$.

Our main theorems establish the analogous result for all Euler pairs. We will prove the theorems both analytically and  combinatorially. We refer to the results in Theorem \ref{T1} as first Beck-type identities and to the result in Theorem \ref{T2} as second Beck-type identity. 

\begin{theorem} \label{T1} If $n,r $ are integers  such that $n \geq 0$ and $r\geq 2$, and $(S_1, S_2)$ is an Euler pair of order $r$, then 
\begin{enumerate}

\item[(i)] $a_r(n)=\ds\frac{1}{r-1}b_r(n)$ \\ \ \\

\item[(ii)] $c_r(n)=\ds\frac{1}{r-1}b_r(n)$. \end{enumerate}\end{theorem}

\begin{theorem} \label{T2} If $n,r $ are integers  such that $n \geq 0$ and $r\geq 2$, and $(S_1, S_2)$ is an Euler pair of order $r$, then    
$c'_r(n)=b'_r(n)$.
\end{theorem}

\begin{example} We continue  with the Euler pair of order $3$ from Example \ref{eg:epair}. 
We have
$$\mathcal{O}_{3}(7) = \{ (7), (5, 1^2), (1^7) \}; \ \mathcal{D}_{3}(7) = \{(7), (5, 1^2), (3^2, 1)\};$$ 
and
$$\mathcal{O}_{1, 3}(7) = \{ (3^2, 1), (3, 1^4) \}; \ \mathcal{D}_{1, 3}(7) = \{(1^7), (3, 1^4)\}.$$

Glaisher's bijection gives us 

\begin{center}
\begin{tabular}{ccc}
$(7)$ & $\stackrel{\varphi_3}{\longrightarrow}$ & $(7)$\\ \ \\ $(5, 1, 1)$ &$\longrightarrow$ & $(5, 1, 1)$ \\ \ \\ $(\underbrace{1, 1, 1}, \underbrace{1, 1, 1}, 1)$ &$\longrightarrow$ & (3, 3, 1)
\end{tabular}
\end{center}

We note that $$a_3(7) = |{O}_{1, 3}(7)| = 2, c_3(7) = |{D}_{1, 3}(7)| = 2, \mbox{ and } b_3(7) = 11 - 7 = 4.$$ Thus, $$\frac{1}{3-1}b_3(7) = a_3(7) = c_3(7).$$

If we restrict to counting different parts in partitions, we see that there are a total of $4$ diferent parts in the partitions of $\mathcal{O}_{3}(7)$ and a total of $5$ different parts in the partitions of $\mathcal{D}_{3}(7)$. Since $\mathcal{T}_3(7) = \{ (3, 1^4) \}$, we have $$b'_3(7) = 5 - 4 = 1 = |\mathcal{T}_3(7)|.$$ 
\end{example}

The analytic proofs of Theorems \ref{T1} and \ref{T2} are similar to the proofs in \cite{A17} and \cite{FT17}, while the combinatorial proofs follow the ideas of \cite{BB19}. However, the generalizations of the proofs in the aforementioned articles to Euler pairs of order $r\geq 2$ are important as establishing the theorems in such generality leads to a multitude of new Beck-type identities. We  reproduce  several Euler pairs listed in \cite{S71}. For each identity $|\OO(n)|=|\D(n)|$ holding for the pair below, there are companion Beck-type identities as in Theorems \ref{T1} and \ref{T2}. 

The following pairs $(S_1, S_2)$ are Euler pairs (of order $2$). 

\begin{enumerate}

\item[(i)] $S_1=\{m\in N : m \not \equiv 0 \pmod 3\}$; 

\noindent $S_2=\{m\in N : m  \equiv 1,5 \pmod 6\}$.

In this case, the identity $|\mathcal{O}_2(n)|=|\mathcal{D}_2(n)|$ is known as Schur's identity.\medskip

\item[(ii)] $S_1=\{m\in N : m  \equiv 2,4,5 \pmod 6\}$; 

\noindent $S_2=\{m\in N : m  \equiv 2,5, 11 \pmod{12}\}$. 

In this case, the identity $|\mathcal{O}_2(n)|=|\mathcal{D}_2(n)|$ is known as G\"ollnitz's identity.\medskip

\item[(iii)] $S_1=\{m\in N : m=x^2+2y^2 \mbox{ for some } x,y\in \Z\}$; 

\noindent $S_2=\{m\in N : m \equiv 1 \pmod 2 \mbox{ and } m=x^2+2y^2 \mbox{ for some } x,y\in \Z\}$.

 \end{enumerate}\medskip

 The following is an Euler pair of order $3$. 
 
 \begin{enumerate}
 
 \item[(iv)] $S_1=\{m\in N : m=x^2+xy+y^2 \mbox{ for some } x,y\in \Z\}$; 

\noindent $S_2=\{m\in N : \gcd(m,3)=1 \mbox{ and } m=x^2+xy+y^2 \mbox{ for some } x,y\in \Z\}$. \end{enumerate}\medskip

The following pairs $(S_1, S_2)$ are Euler pairs of order $r$. 

\begin{enumerate}\medskip

\item[(v)] $S_1=\{m\in N : m  \equiv \pm r \pmod{r(r+1)}\}$; 

\noindent $S_2=\{m\in N : m  \equiv \pm r \pmod{r(r+1)} \mbox{ and } m  \not \equiv \pm  r^2 \pmod{r^2(r+1)} \}$.\medskip

\item[(vi)] $S_1=\{m\in N : m  \equiv \pm r, -1 \pmod{r(r+1)}\}$;

\noindent $S_2=\{m\in N : m \equiv \pm r, -1 \pmod{r(r+1)} \mbox{ and } m  \not \equiv \pm r^2, -r \pmod{r^2(r+1)}\}$.

If $r=2$, this Euler pair becomes G\"ollnitz's pair in (ii) above. \medskip

\item[(vii)] Let $r+1$ be a  prime. 

\noindent $S_1=\{m\in N : m \not \equiv 0 \pmod{r+1}\};$

\noindent $S_2=\{m\in N : m \not \equiv tr, t(r+1) \pmod{r^2+r} \mbox{ for } 1\leq t\leq r \}.$ 

If $r=2$, this Euler pair becomes Schur's pair in (i) above. \medskip

\item[(viii)] Let $p$ be a prime and $r$ a quadratic residue $\pmod p$.

\noindent $S_1=\{m \in \N: \mbox{$m$ quadratic residue} \pmod p\};$

\noindent $S_2=\{m \in \N: m \not \equiv 0 \pmod r \mbox{ and $m$ quadratic residue} \pmod p\}.$

 \end{enumerate}
 Note that each case (v)-(viii) gives infinitely many Euler pairs and therefore leads to infinitely many new Beck-type identities.  We also note that in (vii) we corrected a slight error in (3.4) of  \cite{S71}.
 
 \begin{example} Consider the Euler pair in  (vii) above with $r=4$. 
We have 

$S_1 = \{ m \in N  :  m \not\equiv 0 \ \pmod 5 \};$

\medskip

$S_2=\{m\in N : m \not \equiv 4t, 5t \pmod{20} \mbox{ for } 1\leq t\leq 4 \}.$ \\

 Then $(S_1, S_2)$ is an Euler pair of order $4$ and we have\\

$\mathcal{O}_{4}(7) = \{ (7), (6, 1), (3^2, 1), (3, 2^2), (3, 1^4), (3, 2, 1^2), (2^3, 1), (2^2, 1^3), (2, 1^5), (1^7) \};$

\medskip 

$\mathcal{D}_{4}(7) = \{ (7), (6,1), (3^2, 1), (3, 2^2), (4, 3), (3, 2, 1^2), (2^3, 1), (2^2, 1^3), (4, 2, 1), (4, 1^3)  \}.$

\medskip

We have $\mathcal{O}_{1, 4}(7) = \{ (4, 1^3), (4, 2, 1), (4, 3) \};$
\ \ \  $\mathcal{D}_{1, 4}(7) = \{ (1^7), (2, 1^5), (3, 1^4) \}.$\medskip

We note that $a_4(7) = |{O}_{1, 4}(7)| = 3$, $c_4(7) = |{D}_{1, 4}(7)| = 3$, and $b_4(7) = 40 - 31 = 9$, so $\frac{1}{3}b_4(7) = a_4(7) = c_4(7)$.
\medskip

If we restrict to counting different parts, we see that there are 19 different parts in the partitions of $\mathcal{O}_{4}(7)$ and 21 different parts in the partitions of $\mathcal{D}_{4}(7)$. So $b'_4(7) = 21 - 19 = 2 = |\mathcal{T}_4(7)|$ since $\mathcal{T}_4(7) = \{ (1^7), (2, 1^5) \}$.

\end{example}

\section{Proofs of Theorem \ref{T1} }

\subsection{Analytic Proof}

In this article, whenever we work with $q$-series, we assume that $|q| < 1$. When working with two-variable generating functions, we assume both variables are complex numbers less than $1$ in absolute value. Then all series converge absolutely. The generating functions for $|\D(n)|$ and $|\OO(n)|$ are given by 

\begin{align*} 
\sum_{n = 0}^{\infty}|\D(n)|q^n & = \prod_{a \in S_1} (1 + q^a + q^{2a} + \dotsm +q^{(r-1)a})\\
& = \prod_{a \in S_1} \frac{1 - q^{ra}}{1 - q^a};
\end{align*}
and
\begin{align*}
\sum_{n = 0}^{\infty} |\OO(n)|q^n & = \prod_{{\bf b} \in S_2} \frac{1}{1 - q^{\bf b}}.
\end{align*}

To keep track of the number of parts used, we introduce a second  variable $z$, where $|z| < 1$. Let  $$\D(n; m)=\{\ll\in \D(n) \mid \ll \mbox{ has exactly $m$ parts}\}$$ and $$\OO(n; m)=\{\ll\in \OO(n) \mid \ll \mbox{ has exactly $m$ parts}\}.$$
Then, the generating functions for $|\D(n; m)|$ and $ |\OO(n; m)|$ are
\begin{align*} 
f_{\D}(z,q):=\sum_{n = 0}^{\infty}\sum_{m = 0}^{\infty}|\D(n; m)|z^mq^n & = \prod_{a \in S_1} (1 + zq^a + z^{2}q^{2a} + \dotsm +z^{(r-1)}q^{(r-1)a})\\
& = \prod_{a \in S_1} \frac{1 - z^rq^{ra}}{1 - zq^a};
\end{align*}
and 
\begin{align*}
f_{\OO}(z,q):=\sum_{n = 0}^{\infty}\sum_{m = 0}^{\infty} |\OO(n; m)| z^mq^n & = \prod_{{\bf b} \in S_2} \frac{1}{1 - zq^{\bf b}}.
\end{align*}

To obtain the generating function for  the total number of parts in all partition in $\D(n)$  (respectively $\OO(n)$), 
we take the derivative with respect to $z$ of $f_{\D}(z,q)$ (respectively  $f_{\OO}(z,q)$),  and set $z = 1$. We obtain 


\begin{align*} 
&  \left.\frac{\partial}{\partial z}\right|_{z=1}f_{\D}(z,q)    \\   & \hspace*{1cm}  =  \prod_{a \in S_1} \frac{1 - q^{ra}}{1 - q^a}\sum_{a \in S_1} \frac{-rq^{ra}(1-q^a)+q^a(1-q^{ra})}{(1-q^a)(1-q^{ra})} \\ & \hspace*{1cm}  =  \prod_{a \in S_1} \frac{1 - q^{ra}}{1 - q^a}\sum_{a \in S_1} \left(\frac{q^a}{1-q^a}-\frac{q^{ra}}{1-q^{ra}}-(r-1)\frac{q^{ra}}{1-q^{ra}} \right)\\ & \hspace*{1cm}  =  \prod_{a \in S_1} \frac{1 - q^{ra}}{1 - q^a}\left(\sum_{a \in S_1}  \sum_{\substack{ k=1\\ r\nmid k}}^{\infty}q^{ka}-\sum_{a \in S_1}(r-1)\frac{q^{ra}}{1-q^{ra}} \right)
;
\end{align*}
and

\begin{align*}
\left.\frac{\partial}{\partial z}\right|_{z=1}f_{\OO}(z,q)=\prod_{{\bf b} \in S_2} \frac{1}{1 - q^{\bf b}} \sum_{{\bf b} \in S_2} \frac{q^{\bf b}}{1 - q^{\bf b}}.
\end{align*}

Since $|\D(n)| = |\OO(n)|$, we have 

\begin{align*}
\sum_{n=0}^{\infty} b_r(n) q^n & = \prod_{{\bf b} \in S_2} \frac{1}{1 - q^{\bf b}} \bigg( \sum_{{\bf b} \in S_2} \frac{q^{\bf b}}{1 - q^{\bf b}} -   \sum_{\substack{ a\in S_1\\ k\in \N\\ r\nmid k}}q^{ka}+\sum_{a \in S_1}(r-1)\frac{q^{ra}}{1-q^{ra}}  \bigg).
\end{align*}

Next we see that 

 \begin{align} \label{sets}
 \sum_{{\bf b} \in S_2} \frac{q^{\bf b}}{1 - q^{\bf b}}= \sum_{\substack{ a\in S_1\\ k\in \N\\ r\nmid k}}q^{ka}. 
 \end{align}
 
 We have 
 
 $$\sum_{\substack{ a\in S_1\\ k\in \N\\ r\nmid k}}q^{ka}=\sum_{\substack{ a\in S_1\\ k\in \N}}q^{ka}-\sum_{\substack{ a\in S_1\\ k\in \N}}q^{rka}= \sum_{a\in S_1}\frac{q^a}{1-q^a}- \sum_{a\in S_1}\frac{q^{ra}}{1-q^{ra}}=\sum_{{\bf b}\in S_2}\frac{q^{\bf b}}{1-q^{\bf b}}.$$ The last equality holds because $S_2=S_1\setminus rS_1$. Therefore, \eqref{sets} holds.

Then, the  generating function for $b_r(n)$ becomes

\begin{align*}
\sum_{n=0}^{\infty} b_r(n) q^n &  = \prod_{{\bf b} \in S_2} \frac{1}{1 - q^{\bf b}}  \bigg( (r-1) \sum_{a \in S_1} \frac{q^{ra}}{1 - q^{ra}} \bigg)\\
& = \prod_{a \in S_1} (1 + q^a + q^{2a} + \dotsm +q^{(r-1)a}) \bigg( (r-1) \sum_{a \in S_1} \frac{q^{ra}}{1 - q^{ra}} \bigg).
\end{align*}

Therefore

\begin{align*}
\sum_{n=0}^{\infty} b_r(n) q^n & = \sum_{n = 0}^{\infty} (r-1) |\A(n)|q^n = \sum_{n = 0}^{\infty} (r-1) |\C(n)|q^n.
\end{align*}

\noindent Equating coefficients results in $a_r(n)=\ds\frac{1}{r-1}b_r(n)$ and $c_r(n)=\ds\frac{1}{r-1}b_r(n)$.

\subsection{Combinatorial Proof}
\subsubsection{$b_r(n)$ as the cardinality of a set of marked partitions}\label{b}

We start with another example of Glaisher's bijection.

\begin{example}\label{eg:gl}
We continue  with the Euler pair of order $3$ from Example \ref{eg:epair}, but this time use $n = 11$.\\

$\mathcal{O}_{3}(11) = \{ (11), (7, 1^4), (5^2, 1),  (5, 1^6), (1^{11})   \};$ 

$\mathcal{D}_{3}(11) = \{(11), (9, 1^2),  (7, 3, 1), (5^2, 1), (5, 3^2)  \}.$

\medskip

Thus, $b_3(11)= 27-13=14$. \medskip

Glaisher's  bijection gives us

\begin{center}
\begin{tabular}{ccc}
$(11)$ & $\stackrel{\varphi_3}{\longrightarrow}$ & $(11)$\\ \ \\ $(7, \underbrace{1, 1, 1}, 1)$ &$\longrightarrow$ & $(7, 3, 1)$ \\ \ \\ $(5,5, 1)$ & ${\longrightarrow}$ & $(5,5, 1)$\\ \ \\ $(5, \underbrace{1, 1, 1,} \underbrace{1, 1, 1})$ & ${\longrightarrow}$ & $(5, 3,3)$\\ \ \\ $(\underbrace{\underbrace{1, 1, 1}, \underbrace{1, 1, 1}, \underbrace{1, 1, 1,}} 1, 1)$ & ${\longrightarrow}$ & $(9, 1,1)$
\end{tabular}
\end{center}
\end{example}
From Glaisher's bijection, it is clear that each partition  $\lambda\in \OO(n)$ has at least as many parts as its image $\varphi_r(\lambda)\in \D(n)$.   

When calculating $b_r(n)$,  we sum up the differences in the number of parts in each pair  ($\lambda,\varphi_r(\lambda))$. Write each part $\mu_j$ of  $\mu=\varphi_r(\lambda)$ as  $\mu_j=r^{k_j} \bf a$. Then, $\mu_j$ was obtained by merging $r^{k_j}$ parts equal to $\bf a$ in $\lambda$ and thus contributes an excess of $r^{k_j}-1$ parts to $b_r(n)$. Therefore, the difference between the number of parts of $\lambda$ and the number of parts of $\varphi_r(\lambda)$ is $\ds \sum_{j=1}^{\ell(\varphi_r(\lambda))} (r^{k_j}-1)$.

\medskip

\begin{example} In the setting of  Example  \ref{eg:gl}, we see that $(7,3,1)$ contributes $2$ to $b_3(11)$, $(5,3,3)$ contributes $2+2$ to $b_3(11)$, and $(9,1,1)$ contributes $8$ to $b_3(11)$. Thus, $b_3(11)=2+4+8=14$. \end{example}

\begin{definition} \label{marked} Given $(S_1,S_2)$, an Euler pair of order $r$, we define the set $\MD_{1,r}(n)$ of \textit{marked partitions} of $n$ as the set  of partitions in $\D(n)$ such that exactly one part of the form $r^k\bf a$ with $k\geq 1$  has as index an integer $t$ satisfying $1\leq t\leq r^k-1$. If $\mu\in\D(n)$ has parts $\mu_i=\mu_j=r^k\bf a$, $k\geq 1$, with $i\neq j$, the marked partition in which $\m_i$ has index $t$ is considered different from the marked partition in which $\m_j$ has index $t$.
\end{definition}
Note that marked partitions, by definition, have a non-primitive part. 
Then, from the discussion above we have the following interpretation for $b_r(n)$. 

\begin{proposition} Let $n,r$ be integers such that $n\geq 1$ and $r\geq 2$. Then, $$b_r(n)=|\MD_{1,r}(n)|.$$ \end{proposition}

 \begin{definition} An \textit{$r$-word} $w$ is a sequence of letters from the alphabet $\{0,1, \ldots r-1\}$. The \textit{length} of an $r$-word $w$, denoted $\ell(w)$, is the number of letters in $w$. We refer to \textit{position} $i$ in $w$ as the $i$th entry from the right, where the most right entry is counted as position $0$. 
\end{definition}

Note that leading zeros are allowed and are recorded. For example, if $r=5$, the $5$-words $032$ and $32$ are different even though in base $5$ they both represent $17$. We have $\ell(032)=3$ and $\ell(32)=2$. The empty bit string has length $0$ and is denoted by $\emptyset$.

\begin{definition} Given $(S_1,S_2)$, an Euler pair of order $r$, we define the set  $\MDS(n)$ of \textit{$r$-decorated partitions} as the set of   partitions  in $\D(n)$ with at least one non-primitive part such that exactly one non-primitive part  $r^k\bf a$  is decorated with an index $w$, where $w$ is an $r$-word satisfying   $0\leq\ell(w)\leq k-1$. As in Definition \ref{marked},  if $\mu\in\D(n)$ has  non-primitive parts $\mu_i=\mu_j=r^k\bf a$ with $i\neq j$, the decorated partition in which $\m_i$ has index $w$ is considered different from the decorated partition in which $\m_j$ has index $w$.\end{definition}

Thus, for each  part $\mu_i=r^{k_i}\bf a$ of $\mu\in \MDS(n)$ there are $\ds \frac{r^{k_i}-1}{r-1}$ possible indices and for each partition $\mu\in \MDS(n)$ there are precisely $\ds\frac{1}{r-1} \sum_{j=1}^{\ell(\mu)} (r^{k_j}-1)$ possible decorated partitions with the same parts as $\mu$. \medskip

The discussion above proves the following interpretation for  $\ds\frac{1}{r-1}b_r(n)$. 

\begin{proposition} Let $n,r$ be integers such that $n\geq 1$ and $r\geq 2$. Then, $$\ds\frac{1}{r-1}b_r(n)=|\MDS(n)|.$$ 

\end{proposition}

While it is obvious that $|\MD_{1,r}(n)|=(r-1)|\MDS(n)|$, to see this combinatorially, consider the map $\psi_r:\MD_{1,r}(n) \to \MDS(n)$ defined as follows. If $\ll \in \MD_{1,r}(n)$, then $\psi_r(\ll)$ is the partition in $\MDS(n)$ in which the $r$-decorated part is the same as the marked part in $\ll$. The index of the part of  $\psi_r(\ll)$ is obtained from the index of the part of $\ll$ by writing it in base $r$ and removing the leading digit. 
Clearly, this is a $r-1$ to $1$ mapping. 

\subsubsection{A combinatorial proof for $a_r(n)=\ds\frac{1}{r-1}b_r(n)$}\label{an=bn}

To prove combinatorially   that $a_r(n)=\ds\frac{1}{r-1}b_r(n)$ we establish a one-to-one correspondence between $\A(n)$ and $\MDS(n)$.\medskip

\noindent \textit{From $\MDS(n)$ to $\A(n)$:}\label{part i}

Start with an $r$-decorated partition $\mu \in \MDS(n)$. Suppose the non-primitive part $\mu_i=r^k \bf a$ is decorated with an $r$-word $w$ of length $\ell(w)$. Then, $0\leq\ell(w)\leq k-1$. Let $d_w$ be the decimal value of $w$. We set $d_\emptyset=0$. We  transform $\m$ into a partition $\ll \in \A(n)$ as follows. \medskip

Define $\bar \m$ to be the partition whose parts are all non-primitive parts of $\m$ of the form $\m_j=r^t\bf a$ with $j\leq i$, i.e., all parts $r^t\bf a$ with $t>k$ and, if $\m_i$ is the $p$th part of size $r^k\bf a$ in $\m$, then $\bar \m$ also has $p$ parts equal to $r^k\bf a$. 

Define $\tilde \m$ to be the partition whose parts are all parts of $\m$ that are not in $\bar \m$.

\begin{enumerate}
\item  In $\bar \m$, split one part of size $r^k\bf a$ into $d_w+1$ parts of size $r^{k-\ell(w)}\bf a$ and  $r^k-(d_w+1)r^{k-\ell(w)}$ primitive parts of size $\bf a$. 
 Every other part in $\bar \m$ splits completely into parts of size $r^{k-\ell(w)}\bf a$. Denote the resulting partition by $\bar \ll$.

\item  Let $\tilde \ll=\varphi_r^{-1}(\tilde \m)$.  Thus, $\tilde \ll$ is obtained by splitting all parts in $\tilde \m$  into primitive parts. 
\end{enumerate}

Let $\ll=\bar\ll \cup \tilde \ll$. Then $\ll\in \A(n)$ and its set of non-primitive parts is $\{r^{k-\ell(w)}\bf a\}$. 

\begin{remark} Since $d_w+1\leq r^{\ell(w)}$,  in step 1, the resulting number of primitive parts equal to $\bf a$ is non-negative. Moreover, in $\bar \ll$ there is at least one non-primitive part. 
\end{remark}

\begin{example} We continue  with the Euler pair of order $3$ from Example \ref{eg:epair}. Consider the decorated partition 

\begin{align*}
\mu & =(1215, 135_{02}, 135, 51, 35, 15, 15, 3)\\
& =(3^5 \cdot 5, (3^3 \cdot 5)_{02}, 3^3 \cdot 5, 3\cdot17, 35, 3\cdot 5, 3 \cdot 5, 3 \cdot 1)\in \mathcal{DD}_3(1604).
\end{align*}
We have $k=3, \ell(w)=2, d_w=2$, and \medskip

$\bar{\mu} = (3^5 \cdot 5, 3^3 \cdot 5);$

\medskip

$\tilde{\mu} = (3^3 \cdot 5, 3\cdot17, 35, 3\cdot 5, 3 \cdot 5, 3 \cdot 1).$\\

To create $\bar{\lambda}$ from $\bar{\mu}$:

\begin{enumerate}
\item Part $135=3^3\cdot 5$ splits into three parts of size $15$ and eighteen parts of size $5$. \medskip

\item Part $1215=3^5\cdot 5$ splits into eighty one parts of size $15$. \medskip
\end{enumerate}

This results in $\bar{\lambda} = (15^{84}, 5^{18}).$\\

To create $\tilde{\lambda}$ from $\tilde{\mu}$:\\

All parts in $\tilde{\mu}$ are split into primitive parts. Thus, part $3^3 \cdot 5$ splits into twenty seven parts of size 5, part $3\cdot 17$ splits into three parts of size $17$,  both parts of $3\cdot5$ split into three parts of size $5$ each, and part $3\cdot 1$ splits into three parts of size $1$. Part $35$ is already primitive so remains unchanged.

This results in $\tilde{\lambda} = (35, 17^3, 5^{33}, 1^3)$. Then, setting $\lambda = \bar{\lambda} \cup \tilde{\lambda}$ results in $\lambda=(35, 17^3, 15^{84}, 5^{51}, 1^3)\in\mathcal{O}_{1,3}(1604)$. The non-primitive part is $15=3\cdot 5$. \medskip

\end{example}

\noindent \textit{From $\A(n)$ to $\MDS(n)$:}

Start with a partition $\lambda\in \A(n)$. In $\ll$ there is one and only one non-primitive part  $r^k \bf a$. Let $f$ be the multiplicity of the non-primitive part of $\lambda$. We  transform $\ll$ into an $r$-decorated partition in $\MDS(n)$ as follows. 

Apply Glaisher's bijection to $\ll$ to obtain $\m=\varphi_r(\ll)\in \D(n)$.  Since $\lambda$ has a non-primitive  part, $\mu$ will  have at least one non-primitive part. 

Next, we determine the $r$-decoration of $\mu$. Consider the non-primitive parts $\m_{j_i}$ of $\m$ of the form $r^{t_i} \bf a$  (same $\bf a$ as in the non-primitive part of $\lambda$) and  $t_i\geq k$. Assume $j_1<j_2<\cdots$. For notational convenience, set $\m_{j_0}=0$. Let $h$ be the positive integer such that 
\begin{equation}\label{ineq-h}\sum_{i=0}^{h-1}\m_{j_i}<f\cdot r^k {\bf a} \leq \sum_{i=0}^{h}\m_{j_i}.\end{equation} Then, we will decorate part $\m_{j_h}=r^{t_h}\bf a$.
To determine the decoration, let 
$$N=\frac{\ds\sum_{i=0}^{h-1}\m_{j_i}}{r^k\bf a}.$$ Then, \eqref{ineq-h} becomes $$r^k{\bf a}N<f\cdot r^k{\bf a}\leq r^k{\bf a}N+r^{t_h}{\bf a},$$ which implies $0<f-N\leq r^{t_h-k}$.

Let $d=f-N-1$ and $\ell=t_h-k$. We have $0\leq \ell \leq t_h-1$. Consider the  representation of $d$ in base $r$ and insert leading zeros to form an $r$-word $w$ of length $\ell$. Decorate $\m_{j_h}$ with $w$. The resulting decorated partition is in $\MDS(n)$.

\begin{example} We continue  with the Euler pair of order $3$ from Example \ref{eg:epair}. Consider the partition $\lambda=(35, 17^3, 15^{84}, 5^{51}, 1^3)\in\mathcal{O}_{1,3}(1604)$. The non-primitive part is $15$. We have $k=1, f=84$. 

Glaisher's bijection produces the partition $\mu=(1215, 135^2, 51, 35, 15^2, 3) = (3^5\cdot5, 3^3 \cdot 5, 3^3 \cdot 5, 3 \cdot 17, 35, 3\cdot 5, 3\cdot 5, 3\cdot 1)\in \mathcal D\mathcal D_{3}(1604)$. The parts of the form $3^{r_i}\cdot 5$ with $r_i\geq 1$ are $1215, 135, 135, 15,15$. Since $1215<84(3^1 \cdot 5) \leq 1215 + 135,$ the decorated part will be the first part equal to $135=3^3\cdot 5$. We have $N=1215/15=81$. 

To determine the decoration, let $d=84-81-1=2$ and $\ell=3-1=2$. The base $3$ representation of $d$ is $2$. To form a $3$-word of length $2$, we introduce one leading $0$. Thus, the decoration is $w=02$ and the resulting decorated partition is $(1215, 135_{02}, 135, 51, 35, 15, 15, 3)=(3^5 \cdot 5, (3^3 \cdot 5)_{02}, 3^3 \cdot 5, 3\cdot17, 35, 3\cdot 5, 3 \cdot 5, 3 \cdot 1)\in \mathcal{DD}_{1,3}(1604)$.
\end{example}

\subsubsection{A combinatorial proof for $c_r(n)=\ds\frac{1}{r-1}b_r(n)$}

We note that one can  compose the bijection of section \ref{part i} with the bijection of \cite{FT17} to obtain a combinatorial proof of part (ii) of Theorem \ref{T1}. However, we give an alternative  proof 
 that $c_r(n)=\ds\frac{1}{r-1}b_r(n)$ by establishing a one-to-one correspondence between $\C(n)$ and $\MDS(n)$.  This proof does not involve the bijection of \cite{FT17} and it mirrors the proof of section \ref{an=bn}.\medskip

\noindent \textit{From $\MDS(n)$ to $\C(n)$:}

Start with an $r$-decorated partition $\mu \in \MDS(n)$. Suppose the non-primitive part $\mu_i=r^k \bf a$ is decorated with an $r$-word $w$ of length $\ell(w)$ and decimal value $d_w$. Then, $0\leq\ell(w)\leq k-1$. 
 We transform $\m$ into a partition  $\lambda\in \C(n)$ as follows. \medskip
 
 Let $\bar{\bar \m}$ be the partition whose parts are all non-primitive parts of $\m$ of the form $\m_j=r^t\bf a$ with $j\geq i$, and $k-\ell(w)-1< t\leq k$,
 i.e., all parts $r^t\bf a$ with $k-\ell(w)-1< t<k$ and, if there are $p-1$ parts of size $r^k\bf a$ in $\m$ after the decorated part, then $\bar{\bar \m}$  has $p$ parts equal to $r^k\bf a$. 
 
 Let $\tilde{\tilde \m}$ be the partition whose parts are all parts of $\m$ that are not in $\bar{\bar \m}$. 
 
 In $\bar{\bar \m}$, perform the following steps. 
  \begin{enumerate}
\item  Split one part equal to $r^k\bf a$ into $r(d_w+1)$ parts of size $r^{k-\ell(w)-1}\bf a$ and $m$ primitive parts of size $\bf a$, where $m=r^k-r(d_w+1)r^{k-\ell(w)-1}$.
Apply Glaisher's bijection $\varphi_r$ to the partition consisting of $m$ parts equal to $\bf a$.
 
\item Split all remaining parts of $\bar{\bar \m}$  completely  into parts of size $r^{k-\ell(w)-1}\bf a$. 
\end{enumerate}
 Denote  by $\bar{\bar \ll}$ the partition with parts resulting from steps 1 and 2 above.

Let $\ll=\bar{\bar \ll} \cup \tilde{\tilde \m}$.
Since $r(d_w+1)\geq r$, it follows that $\lambda\in \C(n)$. The  part repeated at least $r$ times is $r^{k-\ell(w)-1}\bf a$. 

\medskip

\begin{remark} (i) Since $d_w+1\leq r^{\ell(w)}$, the splitting in step 1 can be performed.  

(ii) Note that $r \mid m=r^{k-\ell(w)}(r^{\ell(w)}-(d_w+1))$. Thus, if $w\neq \emptyset$, after applying Glaisher's bijection $\varphi_r$ to the partition consisting of $m$ parts equal to $\bf a$, we obtain parts $r^j\bf a$ with $k-\ell(w)\leq j<k$. Since in $\tilde{\tilde \mu}$, all parts of the form $r^i\bf a$ have $i \geq k$ or $i\leq k-\ell(w)-1$, Glaisher's bijection in step (1) does not create parts equal to any parts in $\tilde{\tilde \mu}$. 

(iii) If $w=\emptyset$, then, in step (1), we have $m=0$ and  $r^k\bf a$ splits into $r$ parts equal to $r^{k-1}\bf a$.  
\end{remark}

\begin{example} We continue  with the Euler pair of order $3$ from Example \ref{eg:epair}. Consider the partition $\mu = (32805, (10935)_{0120}, 10935, 1215, 45, 45, 25, 9, 3) = (3^8 \cdot 5, (3^7\cdot 5)_{0120}, 3^7 \cdot 5, 3^5 \cdot 5, 3^2 \cdot 5, 3^2 \cdot 5, 25, 3^2 \cdot 1, 3 \cdot 1) \in \mathcal{DD}_3(56017)$. Then the decorated part is $\mu_2 = 3^7 \cdot 5$ and the decoration is $w = 0120$. We have $k=7$, $\ell(w) = 4$, $d_w = 15$. So\\

$\bar{\bar{\mu}} = (3^7 \cdot 5, 3^7 \cdot 5, 3^5 \cdot 5);$

\medskip

$\tilde{\tilde{\mu}} = (3^8 \cdot 5, 3^2 \cdot 5, 3^2 \cdot 5, 25, 3^2 \cdot 1, 3 \cdot 1).$

\begin{enumerate}

\item $3^7 \cdot 5$ splits into 

\begin{itemize}
\item $r(d_w + 1) = 48$ parts of $3^2 \cdot 5$ and 
\item $m = r^k - r(d_w + 1)r^{k- \ell(w) - 1} = 3^7 - 48(3^2) = 1755$ parts of $5$.
\end{itemize}

The 1755 parts of 5 merge into two parts of 3645, one part of 1215, and two parts of 135.

\item $3^7 \cdot 5$ splits into two hundred and forty three parts of $3^2 \cdot 5$ and $3^5 \cdot 5$ splits into twenty seven parts of $3^2 \cdot 5$.

\end{enumerate}

This results in\\

$\bar{\bar{\lambda}} = (3645^2, 1215, 135^2, 45^{318});$ 

\medskip

$\lambda = \bar{\bar{\lambda}} \cup \tilde{\tilde{\mu}} = (32805, 3645^2, 1215,135^2, 45^{320}, 25, 9, 3) \in \mathcal D_{1,3}(56017)$. The part repeated at least three times is $45=3^2\cdot 5$.

\end{example}\medskip

\noindent \textit{From $\C(n)$ to $\MDS(n)$:}

Start with a partition $\lambda\in \C(n)$. Then, among the parts of $\lambda$, there is one and only one  part  that is repeated at least $r$ times. Suppose this  part is  $r^k \bf a$, $k\geq 0$,  and denote by  $f\geq r$ its  multiplicity in $\lambda$. As in  Glaisher's bijection, we merge repeatedly parts of $\ll$ that are repeated at least $r$ times to obtain  $\m\in \D(n)$. Since $\lambda$ has a part repeated at least $r$ times, $\mu$ will  have at least one  non-primitive  part.

Next, we determine the decoration of $\mu$.  In this case, we want to work with the parts of $\mu$ from the right to the left (i.e., from smallest to largest part). Let $\tilde{\mu}_q=\mu_{\ell(\mu)-q+1}$. 
 Consider the parts $\tilde\m_{j_i}$ of the form $r^{t_i} \bf a$ with   $t_i\geq k$. If $t_1\leq t_2\leq\cdots$, we have $j_1<j_2<\cdots$. 

As before, we set $\tilde{\m}_{j_0}=0$. Let $h$ be the positive integer such that 
\begin{equation}\label{ineq1-h}\sum_{i=0}^{h-1}\tilde\m_{j_i}<f\cdot r^k {\bf a} \leq \sum_{i=0}^{h}\tilde\m_{j_i}.\end{equation} Then, we will decorate part $\tilde\m_{j_h}=r^{t_h}\bf a$. 
To determine the decoration, let   \begin{equation} \label{Nh} N=\frac{\ds\sum_{i=0}^{h-1}\tilde\m_{j_i}}{r^k\bf a}.\end{equation} Then, \eqref{ineq1-h} becomes $$r^k{\bf a}N<f\cdot r^k{\bf a}\leq r^k{\bf a}N+r^{t_h}{\bf a},$$ which implies $0<f-N\leq r^{t_h-k}$.

Let $\ds d=\frac{f-N}{r}-1$ and $\ell=t_h-k-1$. Note that, by construction, $t_h>k$, and therefore   $0\leq \ell \leq t_h-1$. Consider the  representation of $d$ in base $r$ and insert leading zeros to form an $r$-word $w$ of length $\ell$. Decorate $\tilde\m_{j_h}$ with $w$. The resulting decorated partition (with parts written in non-increasing order) is in $\MDS(n)$. \medskip

\begin{remark} To see that $f-N$ above is always divisible by $r$, note that if $f=qr+t$ with $q,t \in \Z$ and $0\leq t<r$, then there are $t$ terms equal to $r^k\bf a$ in the numerator of $N$. All other terms, if any,  are divisible by $r^{k+1}\bf a$. Therefore, the remainder of $N$ upon division by $r$ is $t$.  
\end{remark} 

\begin{example} We continue  with the Euler pair of order $3$ from Example \ref{eg:epair}. Consider the partition $\lambda = (32805, 3645^2, 1215, 135^2, 45^{320}, 25, 9, 3) \in \mathcal{D}_{1,3}(56017)$. The part repeated at least three times is $45=3^2\cdot 5$. We have $k=2$ and $f = 320$. 

Applying Glaisher's bijection to $\lambda$ results in\\

$\mu = \varphi_{3}(\lambda) = (3^8 \cdot 5, 3^7\cdot 5, 3^7 \cdot 5, 3^5 \cdot 5, 3^2 \cdot 5, 3^2 \cdot 5, 25, 3^2, 3 \cdot 1) \in \mathcal{D}_{3}(56017)$.\\

The parts of the form $3^{t_i}\cdot 5$ with $t_i \geq 2$ are $3^2\cdot 5, 3^2\cdot 5, 3^5\cdot 5, 3^7\cdot 5, 3^7\cdot 5, 3^8\cdot 5$. Since $3^2\cdot 5 + 3^2\cdot 5 + 3^5\cdot 5 + 3^7\cdot 5 < 320 \cdot 3^2 \cdot 5 \leq 3^2\cdot 5 + 3^2\cdot 5 + 3^5\cdot 5 + 3^7\cdot 5 + 3^7\cdot 5 $, the decorated part will be the second part (counting from the right) equal to $3^7 \cdot 5 = 10935$. We have  $N = \ds\frac{3^2\cdot 5 + 3^2\cdot 5 + 3^5\cdot 5 + 3^7\cdot 5}{3^2 \cdot 5} = 272$. Thus $d = \ds\frac{320-272}{3} - 1 = 15$ and $\ell = 7 - 2 - 1 = 4$. The base $3$ representation of $d$ is $120$. To form a $3$-word of length $4,$ we introduce one leading $0$. Thus, the decoration is $w=0120$ and the resulting decorated partition is 

\begin{align*}
\mu = & (32805, (10935)_{0120}, 10935, 1215, 45, 45, 25, 9, 3)\\ 
= & (3^8 \cdot 5, (3^7\cdot 5)_{0120}, 3^7 \cdot 5, 3^5 \cdot 5, 3^2 \cdot 5, 3^2 \cdot 5, 25, 3^2 \cdot 1, 3 \cdot 1) \in \mathcal{DD}_3(56017).
\end{align*}

\end{example}

\section{Proofs of Theorem \ref{T2} }

\subsection{Analytic Proof}

We create a bivariate generating function to keep track of the number of different parts in partitions in $\OO(n)$, respectively $\D(n)$.

We denote by $\mathcal O'_r(n;m)$ the set of partitions of $n$ with parts from $S_2$ using $m$ different parts. We denote by $\mathcal D'_r(n;m)$ the set of partitions of $n$ with parts from $S_1$ using $m$ different parts and allowing parts to repeat no more than $r-1$ times. Then,

\begin{align*}
f_{\mathcal O'_r}(z,q):=\displaystyle \sum_{n = 0}^{\infty} \sum_{m=0}^\infty| \mathcal O'_r(n;m)|z^mq^n & = \prod_{{\bf b} \in S_2} (1 + zq^{\bf b} + zq^{2\bf b} + \dotsm)\\
& = \prod_{{\bf b}  \in S_2} \bigg( 1 + \frac{zq^{\bf b} }{1- q^{\bf b} } \bigg),
\end{align*}
and \begin{align*}
f_{\mathcal D'_r}(z,q):=\displaystyle \sum_{n = 0}^{\infty}\sum_{m=0}^\infty |\mathcal D'_r(n;m)|z^mq^n & = \prod_{a \in S_1} (1 + zq^a + \dotsm +zq^{(r-1)a})\\ & = \prod_{a \in S_1} \bigg( 1 + \frac{zq^a-zq^{ra}}{1- q^a } \bigg).
\end{align*}

To obtain the generating function for  the total number of different parts in all partitions in $\OO(n)$  (respectively $\D(n)$), 
we take the derivative with respect to $z$ of $f_{\mathcal O'_r}(z,q)$ (respectively  $f_{\mathcal D'_r}(z,q)$),  and set $z = 1$. We obtain 

\begin{align*} 
 \left.\frac{\partial}{\partial z}\right|_{z=1}f_{\mathcal O'_r}(z,q) &=\sum_{{\bf b}  \in S_2} \frac{q^{\bf b} }{1 - q^{\bf b} } \prod_{{\bf c} \in S_2, {\bf c} \neq {\bf b} } \bigg(1 + \frac{q^{\bf c}}{1 - q^{\bf c}} \bigg)\\ &=\sum_{{\bf b}  \in S_2} \frac{q^{\bf b} }{1 - q^{\bf b} } \prod_{{\bf c} \in S_2, {\bf c} \neq {\bf b} } \bigg( \frac{1}{1 - q^{\bf c}} \bigg)\\
& = \prod_{{\bf b}  \in S_2} \frac{1}{1 - q^{\bf b} } \sum_{{\bf b}  \in S_2} q^{\bf b},
\end{align*}

and

\begin{align*} 
 \left.\frac{\partial}{\partial z}\right|_{z=1}f_{\mathcal D'_r}(z,q) &= \sum_{a \in S_1} \frac{q^a-q^{ra}}{1-q^a} \prod_{d \in S_1, d\neq a} \bigg(1 + \frac{q^d-q^{rd}}{1-q^d}\bigg)\\ & = \sum_{a \in S_1} \frac{q^a-q^{ra}}{1-q^a} \prod_{d \in S_1, d\neq a}\frac{1-q^{rd}}{1-q^d}\\ 
& =\prod_{a \in S_1}\frac{1-q^{ra}}{1-q^a} \sum_{a \in S_1} \frac{q^a - q^{ra}}{1 - q^{ra}}.
\end{align*}

Since $|\D(n)| = |\OO(n)|$, we have 

\begin{align*} 
 \sum_{n = 0}^{\infty} b'_r(n) q^n 
= \prod_{a \in S_1}\frac{1-q^{ra}}{1-q^a}\bigg( \sum_{a \in S_1} \frac{q^a}{1 - q^{ra}} -\sum_{a \in S_1} \frac{q^{ra}}{1 - q^{ra}} - \sum_{{\bf b}  \in S_2} q^{\bf b}  \bigg).
\end{align*}

Moreover, 
\begin{align*}
 \sum_{a \in S_1} \frac{q^a}{1 - q^{ra}} -\sum_{a \in S_1} \frac{q^{ra}}{1 - q^{ra}}  & = \bigg( \sum_{a \in S_1} q^a +  \sum_{a \in S_1} \frac{q^{(r+1)a}}{1 - q^{ra}} \bigg) - \bigg( \sum_{a \in rS_1} q^a + \sum_{a \in S_1} \frac{q^{2ra}}{1 - q^{ra}} \bigg)  \\
 & = \bigg( \sum_{a \in S_1} q^a - \sum_{a \in rS_1} q^a \bigg) + \bigg(\sum_{a \in S_1} \frac{q^{(r+1)a}}{1 - q^{ra}} - \sum_{a \in S_1} \frac{q^{2ra}}{1 - q^{ra}}\bigg)\\
 & = \sum_{{\bf b}  \in S_2} q^{\bf b} + \sum_{a \in S_1} \frac{q^{(r + 1)a} - q^{2ra}}{1 - q^{ra}},
 \end{align*}
 the last equality occurring because $S_1 = S_2 \sqcup rS_1$. 
 
 Therefore, 
\begin{align*} 
 \sum_{n = 0}^{\infty} b'_r(n) q^n & =  \prod_{a \in S_1}\frac{1-q^{ra}}{1-q^a}\sum_{a \in S_1} \frac{q^{(r + 1)a} - q^{2ra}}{1 - q^{ra}}\\ & = \sum_{a \in S_1} \frac{q^{(r + 1)a} + q^{(r + 2)a} + \dotsm +q^{(2r-1)a}}{1 + q^a + \dotsm + q^{(r-1)a}} \prod_{a \in S_1} (1 + q^a + \dotsm + q^{(r-1)a})\\ & = \sum_{a \in S_1} (q^{(r + 1)a} + q^{(r + 2)a} + \dotsm +q^{(2r-1)a}) \prod_{d \in S_1, c\neq a} (1 + q^d + \dotsm + q^{(r-1)d})\\
& = \sum_{n = 0}^{\infty} c'_r(n) q^n.
 \end{align*}

\medskip

\subsection{Combinatorial Proof}

\subsubsection{$b'_r(n)$ as the cardinality of a set of overpartitions} \label{b1}

As in section \ref{b}, we use Glaisher's bijection and calculate $b'_r(n)$ by summing up the difference between the number of different parts of $\varphi_r(\lambda)$ and the number of different parts of $\lambda$ for each partition $\lambda\in \OO(n)$. For a given ${\bf a} \in S_2$, each part in $\varphi_r(\ll)$ of the form $r^k\bf a$, $k\geq 0$, is obtained from $\ll$ by merging $r^k$ parts equal to $\bf a$. Therefore, the contribution to $b'_r(n)$ of each $\m\in \D(n)$ equals $$\sum_{\substack{{\bf a}\in S_2 \\ {\bf a} \mbox{ \small{part of} } \varphi_r^{-1}(\m)}}(m_{\m}({\bf a})-1),$$ where 
$$m_{\m}({\bf a})=|\{t\geq 0 \mid r^t{\bf a} \mbox{ is a part of } \m \}|. $$

Next, we define a set of overpartitions. An overpartition is a partition in which the last appearance of a part may be overlined. For example,  $(5, \bar 5, 3,3, \bar 2, 1, 1, \bar 1)$ is an overpartition of $21$. We denote by  $\OD(n)$  the set of overpartitions of $n$ with  parts in $S_1$ repeated at most $r-1$ times in which \textit{exactly one} part is overlined and such that  part $r^s\bf a$ with $s \geq 0$  may be overlined only if there is a part $r^t\bf a$ with $t<s$. In particular, no primitive part can be overlined. Note that when we count parts in an overpartition, the overlined part contributes to the multiplicity.  The discussion above proves the following interpretation of $b'_r(n)$.

\begin{proposition} Let $n\geq 1$. Then, $b'_r(n)=|\OD(n)|.$ \end{proposition}

\subsubsection{A combinatorial proof for $c'_r(n)=b'_r(n)$} We establish a one-to-one correspondence between $\OD(n)$ and $\T(n)$.\medskip

\noindent \textit{From  $\OD(n)$ to $\T(n)$:} 

Start with an overpartition $\mu \in \OD(n)$. Suppose the overlined part is $\mu_i=r^s\bf a$. Then there is a part $\mu_p=r^t\bf a$ of $\mu$ with $t<s$. Let $k$ be the largest non-negative integer such that $r^k\bf a$ is a part of $\mu$ and $k<s$. To obtain $\lambda \in \T(n)$ from $\mu$, split $\mu_i$ into $r$ parts equal to $r^k\bf a$ and $r-1$ parts equal to $r^j\bf a$  for each $j=k+1, k+2, \ldots, s-1$.

\begin{example} We continue  with the Euler pair of order $3$ from Example \ref{eg:epair}. Let\\

$\mu = (3^8 \cdot 5, 3^7 \cdot 5, \overline{3^7 \cdot 5}, 3^5 \cdot 5, 3^2 \cdot 5, 3^2 \cdot 5, 25, 3^2 \cdot 1, 3 \cdot 1) \in \ODT(56017)$.\\

Then $k=5$ and $3^7 \cdot 5$ splits into three parts equal to $3^5\cdot 5$ and two parts equal to $3^6 \cdot 5$. Thus, we obtain the partition
\begin{align*}
\lambda & = (3^8 \cdot 5, 3^7 \cdot 5, 3^6 \cdot 5, 3^6 \cdot 5, 3^5 \cdot 5, 3^5 \cdot 5, 3^5 \cdot 5, 3^5 \cdot 5, 3^2 \cdot 5, 3^2 \cdot 5, 25, 3^2 \cdot 1, 3 \cdot  1)\\
& \in \mathcal{T}_3(56017).
\end{align*}
The part repeated more than three times but less than six times is $3^5\cdot 5$.

\end{example}

\noindent \textit{From $\T(n)$ to $\OD(n)$:}

Start with a partition $\lambda \in \T(n)$. Suppose $r^k\bf a$ is the part repeated more than $r$ times but less than $2r$ times.   Let $\m=\varphi_r(\ll)\in \D(n)$. Overline the smallest part of $\mu$ of form $r^t \bf a$ with $t>k$.  The resulting overpartition is  in $\OD(n)$.

\begin{example} We continue  with the Euler pair of order $3$ from Example \ref{eg:epair}. Let
\begin{align*}
\lambda & = (3^8 \cdot 5,  3^7 \cdot 5,  3^6 \cdot 5, 3^6 \cdot 5, 3^5 \cdot 5, 3^5 \cdot 5, 3^5 \cdot 5, 3^5 \cdot 5, 3^2 \cdot 5, 3^2 \cdot 5, 25, 3^2 \cdot 1, 3 \cdot  1)\\
& \in \mathcal{T}_3(56017).
\end{align*}The part repeated more than three times but less than six times is $3^5\cdot 5$.
We have $k = 5$. Merging by Glaisher's bijection, we obtain\\

$\mu = (3^8 \cdot 5, 3^7 \cdot 5, 3^7 \cdot 5, 3^5 \cdot 5, 3^2 \cdot 5, 3^2 \cdot 5, 25, 3^2 \cdot 1, 3 \cdot 1) \in \mathcal{D}_3(56017).$\\

The smallest part of $\mu$ of the form $r^t\bf{a}$ with $t > k = 5$ is $3^7\cdot 5$. Thus we obtain the overpartition\\

$\mu = (3^8 \cdot 5, 3^7 \cdot 5, \overline{3^7 \cdot 5}, 3^5 \cdot 5, 3^2 \cdot 5, 3^2 \cdot 5, 25, 3^2 \cdot 1, 3 \cdot 1) \in \ODT(56017).$

\end{example}

\begin{remark} \label{rem}  We could have obtained the  transformation above from the combinatorial proof of part (ii) of Theorem \ref{T1}. In the transformation from $\C(n)$ to $\MDS(n)$, if part $r^k\bf a$ is the part repeated more than $r$ times but less than $2r$ times, we have $f=r+s$ for some $1\leq s\leq r-1$, $h=s+1$, and $N=s$. Thus $d=0$ and the decorated part is the last occurrence of smallest part in the transformed partition $\mu$ that is of the form $r^t\bf a$ with $t>k$.  Thus, in $\mu$, the decorated part $r^{t} \bf a$ is decorated with an $r$-word consisting of all zeros and of length $t-k-1$, one less than the difference in exponents of $r$ of the decorated part and the next smallest part with the same $\bf a$ factor. Since in this case the decoration of a partition in $\MDS(n)$ is completely determined by the part being decorated, we can simply just overline the part.  \end{remark}

\section{Concluding remarks}

In this article we proved first and second Beck-type identities for all Euler pairs $(S_1,S_2)$ of order $r\geq 2$. Euler pairs of order $r$ satisfy $rS_1\subseteq S_1$ and $S_2=S_1\setminus rS_1$. 
Thus, we established Beck-type identities accompanying all partition identities of the type given in Theorem \ref{sub-thm}. 

At the end of \cite{S71}, Subbarao mentions that the characterization of Euler pairs of order $r$ given by Theorem \ref{sub-thm} can be extended to vector partitions.  The corrected statement for partitions of multipartite numbers is given in  \cite[Theorem 12.2]{A98} and indeed it has Beck-type companion identities as we explain below. 

A multipartite (or $s$-partite) number $\frak n=(n_1, n_2, \ldots, n_s)$ is an $s$ tuple of non-negative integers, not all $0$. We view multipartite numbers as vectors and refer to $n_1, n_2, \ldots, n_s$ as the entries of $\frak n$. 

A multipartition (or vector partition) $\xi=(\xi^{(1)}, \xi^{(2)},\ldots, \xi^{(t)})$ of $\frak n$ is a sequence of multipartite numbers in non-increasing lexicographic order satisfying   $$\frak n=\xi^{(1)}+ \xi^{(2)}+\cdots +\xi^{(t)}.$$ We refer to $\xi^{(i)}$, $1\leq i\leq t$ as the multiparts (or vector parts) of the multipartition $\xi$ and to the number of multiparts $t$ of $\xi$ as the length of $\xi$ which we denote by $\ell(\xi)$. 

Let $S_1$ and $S_2$ be sets of positive integers. Given a multipartition $\xi$ of $\frak n$ with all entries of all multiparts in $S_1$, we say that a multipart $\xi^{(i)}$ of $\xi$ is \textit{primitive} if at least one entry of $\xi^{(i)}$ is in $S_2$. Otherwise, the multipart is called \textit{non-prmitive}. We denote by  $\mathcal{VD}_r(\frak n)$ the set of multipartitions $\xi=(\xi^{(1)}, \xi^{(2)},\ldots, \xi^{(t)})$ of $\frak n=(n_1, n_2, \ldots, n_s)$ with all entries of all multiparts in $S_1$ and such that all multiparts are repeated at most $r-1$ times. We denote by  $\mathcal{VO}_r(\frak n)$ the set of multipartitions $\eta=(\eta^{(1)}, \eta^{(2)}, \ldots, \eta^{(u)})$ of $\frak n=(n_1, n_2, \ldots, n_s)$ such that each multipart $\eta^{(i)}$ of $\eta$ is primitive.

Then, Andrews \cite{A98} gives the following theorem mentioning that its proof can be constructed similar to the proof using ideals of order $1$ for the analogous result for regular partitions. 

\begin{theorem} \label{SA} Let $S_1$ and $S_2$ be sets of positive integers. Then $$|\mathcal{VD}_r(\frak n)|=|\mathcal{VO}_r(\frak n)|$$ if and only if $(S_1,S_2)$ is an Euler pair of order $r$, i.e., $rS_1\subseteq S_1$ and $S_2=S_1-rS_1$.

\end{theorem}

We note that Glaisher's bijection can be extended to prove Theorem \ref{SA} combinatorially. The Glaisher type transformation, $v\varphi_r$, from $\mathcal{VO}_r(\frak n)$ to $\mathcal{VD}_r(\frak n)$ repeatedly merges $r$ equal multiparts (as addition of vectors) until there are no multiparts repeated more than $r-1$ times. The transformation from $\mathcal{VD}_r(\frak n)$  to $\mathcal{VO}_r(\frak n)$ takes each non-primitive multipart  (all its entries are from $rS_1$) and splits it into $r$ equal multiparts, repeating the process until the obtained multiparts are primitive. The remark below is the key to adapting the combinatorial proofs of Theorems \ref{T1} and \ref{T2} to proofs of Beck-type identities for multipartitions. 

\begin{remark} Let $\eta \in \mathcal{VO}_r(\frak n)$ and $\xi=v\varphi_r(\eta)\in \mathcal{VD}_r(\frak n)$. Then multipart $\xi^{(i)}$ of $\xi$ was obtained by merging $r^k$ multiparts of $\eta$ if and only if, when writing all entries of $\xi^{(i)}$ in the form $r^j{\bf a}$, the smallest exponent of $r$ in all entries of $\xi^{(i)}$ is $k$. 
\end{remark}

To formulate Beck-type identities for multipartitions, let  $vb_r(\frak n)$ be the difference between the number of multiparts in all multipartitions in $\mathcal{VO}_r(\frak n)$ and the number of multiparts in all multipartitions in $\mathcal{VD}_r(\frak n)$. Similarly,  let  $vb'_r(\frak n)$ be  the difference in the total number of different multiparts in all multipartitions in $\mathcal{VD}_r(\frak n)$  and the total number of different multiparts in all multipartitions in $\mathcal{VO}_r(\frak n)$. Then, we have the following Beck-type identities for multipartions.

\begin{theorem} \label{last} Suppose $(S_1,S_2)$ is an Euler pair of order $r\geq 2$  and let $\frak n$ be a multipartite number. Then 

\begin{enumerate}

\item[(i)] $\ds \frac{1}{r-1}vb_r(\frak n)$ equals the number of multipartitions $\xi=(\xi^{(1)}, \xi^{(2)},\ldots, \xi^{(t)})$ of $\frak n$ with all entries of all multiparts in $S_1$ and such that exactly one multipart is repeated at least $r$ times. Moreover, $\ds \frac{1}{r-1}vb_r(\frak n)$ equals the number of multipartitions $\eta=(\eta^{(1)}, \eta^{(2)}, \ldots, \eta^{(u)})$ of $\frak n$ such that each multipart $\eta^{(i)}$ of $\eta$ has all entries in $S_1$ and only one multipart, possibly repeated, is non-primitive. \medskip

\item[(ii)]  $vb'_r(\frak n)$ equals   the number of multipartitions $\xi=(\xi^{(1)}, \xi^{(2)},\ldots, \xi^{(t)})$ of $\frak n$ with all entries of all multiparts in $S_1$ and such that exactly one multipart is  repeated more than $r$ times but less than $2r$ times. \end{enumerate}\end{theorem}

The combinatorial proofs of these statements follow the combinatorial proofs of Theorems \ref{T1} and \ref{T2} with  all references to expression of the form ``part $r^k{\bf a}$ in the partition $\mu$" changed to ``multipart of multipartition $\mu$ in which $r^k{\bf a}$ is the entry with the smallest exponent of $r$ (among the entries of the multipart)."

To our knowledge, there is no  analytic proof of Theorem \ref{last}. 

\section*{Acknowledgements}
We are grateful to the anonymous referees for suggestions that improved the exposition of the article. In particular, one referee suggested the short proof of \eqref{sets}, and another referee alerted us to the correct statement of Subbarao's theorem for vector partitions.

\bigskip


\end{document}